\documentclass{au}
\usepackage{amsmath}
\newtheorem{theorem}{Theorem}
\newtheorem{lemma}{Lemma}

\newtheorem{observation}{Observation}

\newtheorem{definition}{Definition}

\def\proof{\noindent{\bf Proof. }}

\def\marcafdem{\vrule width 1.2ex height 1.1ex depth 0.1ex}
\def\qed{\hbox{}\nobreak\hfill\hbox{\marcafdem}\par}

\def\endproof{\qed}

\long\def\BEGINCOMMENT #1\ENDCOMMENT{\relax}

\newcommand{\maps}\longrightarrow
\newcommand{\cmaps}\Longrightarrow
\newcommand{\vocab}{\tau}

\newcommand{\best}{\operatorname{{\bf B}}}
\newcommand{\aest}{\operatorname{{\bf A}}}
\newcommand{\csp}{\operatorname{CSP}}

\newcommand{\Inv}{\operatorname{Inv}}
\newcommand{\Pol}{\operatorname{Pol}}

\newcommand{\pr}{\operatorname{pr}}

\newcommand{\sg}{\operatorname{Sg}}

\newcommand{\dom}{\operatorname{dom}}

\newcommand{\algB}{\operatorname{{\mathcal B}}}
\newcommand{\algA}{\operatorname{{\mathcal A}}}
\newcommand{\algR}{\operatorname{{\mathcal R}}}
\newcommand{\algS}{\operatorname{{\mathcal S}}}
\newcommand{\algD}{\operatorname{{\mathcal D}}}
\newcommand{\algC}{\operatorname{{\mathcal C}}}
\newcommand{\algE}{\operatorname{{\mathcal E}}}
\newcommand{\algH}{\operatorname{{\mathcal H}}}

\begin{document}

\title{CD(4) has bounded width}

\author{Catarina Carvalho} 
\author{V\'{\i}ctor Dalmau}
\author{Petar Markovi\'c}
\author{Mikl\'os Mar\'oti}
\address{Department of Computer Science, Durham University, United Kingdom}
\email{c.a.carvalho@durham.ac.uk}
\address{Departament de Tecnologia, Universitat Pompeu Fabra, Barcelona, Spain}
\email{victor.dalmau@upf.edu}
\address{Department of Mathematics and Informatics, University of Novi Sad, Serbia}
\email{pera@im.ns.ac.yu}
\address{Bolyai Institute, University of Szeged, Hungary}
\email{mmaroti@math.u-szeged.hu}
\date{September, 2007}
\thanks{The third author was supported by the grant no. 144011G of the Ministry of Science and Environment of Serbia.}

\begin{abstract}
We prove that the constraint languages invariant under a short sequence of J\'onsson terms (containing at most three non-trivial ternary terms) are tractable by showing that they have bounded width. This improves the previous result by Kiss and Valeriote \cite{Kiss/Valeriote} and presents some evidence that the Larose-Zadori conjecture \cite{Larose/Zadori} holds in the congruence-distributive case.
\end{abstract}

\maketitle

\section{Introduction}

In recent years, universal algebra has proven to be very useful in
the study of the computational complexity of the constraint
satisfaction problem. For every relational structure $\best$, the
constraint satisfaction problem (CSP) associated to $\best$,
$\csp(\best)$, is the following computational problem: given a
finite structure $\aest$, determine whether $\aest$ is homomorphic
to $\best$. Many computational problems, coming from areas as
diverse as artificial intelligence, scheduling, graph theory,
database theory, and others can be formulated, in a natural way, as
a constraint satisfaction problem. From a computational complexity
point of view the importance of the CSP was first pointed out by
Feder and Vardi \cite{Feder/Vardi} who have shown that if the class of
constraint satisfaction problems, in its logic formulation, is
slightly generalized in several different ways then we obtain a
class of problems which is essentially as rich as the whole of NP.
This fact motivates the dichotomy question "are there CSPs that are
not solvable in polynomial time nor NP-complete?" which despite
considerable effort still remains open.

The groundbreaking work of Jeavons, Cohen, and Gyssens
\cite{Jeavonsetal} successively developed and refined by Bulatov,
Jeavons, and Krokhin \cite{Bulatov/Jeavons/Krokhin}
and Larose and Tesson \cite{Larose/Tesson} has shown strong ties between CSP and
universal algebra. In particular, it has been shown that the
computational complexity of $\csp(\best)$ is uniquely determined by
the algebra ${\mathcal A}_{\best}$ which has the same universe as
$\best$ and whose basic operations are the polymorphisms of the
relations in $\best$. A good deal of recent results on the
complexity of the CSP are due to this link (see the survey of
Bulatov, Jeavons and Krokhin \cite{Bulatov/Jeavons/Krokhin:survey} for an overview). It is worth mentioning that all this activity has spurred development of universal algebra itself, witnessed mostly in the development of all sorts of new Mal'cev-style conditions. Some examples of results of this sort are \cite{mar-mck}, \cite{BIMMVW} and \cite{Barto/Kozik/Niven}.

There are basically two algorithms, or rather algorithmic
principles, for CSPs. The first one is linked to the few subalgebras
property studied in \cite{Dalmau} and \cite{IMMVW}.

The second one, central to this paper, is called the $k$-consistency
algorithm and gives rise to the notion of bounded width (see the
recent survey \cite{Bulatov/Krokhin/Larose} by Bulatov, Krokhin, and Larose for a nice overview
of bounded width). In a nutshell, for every fixed $k>0$, the
$k$-consistency algorithm is an iterative, polynomial-time,
algorithm that computes a set $H$ of partial homomorphisms from
$\aest$ to $\best$ satisfying the condition that every complete
homomorphism (if it exists) from $\aest$ to $\best$ must have all its
$k$-ary projections in $H$ (see section \ref{sect:prelim} for precise
definitions). When the set $H$ returned by the $k$-consistency
algorithm is empty we have a guarantee that there is no homomorphism
from $\aest$ to $\best$. Those relational structures $\best$ for
which there exists some $k>0$ such that the converse also holds are
said to have bounded width. Consequently, if $\best$ has bounded
width then it is possible to use the $k$-consistency algorithm for
some $k>0$ to solve correctly $\csp(\best)$ in polynomial time. A
most important question in the area is to determine which structures
$\best$ has bounded width, which is equivalent to delineate the
reach of the $k$-consistency algorithm as a tool to solve CSPs. In
this study, universal algebra has played a major role. Larose and
Z\'adori \cite{Larose/Zadori} have shown that if $\best$ has bounded width then its
associated algebra ${\mathcal A}_{\best}$ generates a variety that
omits Hobby-McKenzie types ${\mathbf 1}$ and ${\mathbf 2}$ (also Bulatov has proved an essentially
equivalent statement in \cite{Bulatovrelational}). It has  been conjectured in \cite{Larose/Zadori} that
this condition is also sufficient. The Larose-Zadori conjecture
would imply, in particular, that any structure $\best$ whose
associated algebra ${\mathcal A}_{\best}$ is in CD has bounded
width. This has only been verified for algebras containing a
near-unanimity term by Feder and Vardi \cite{Feder/Vardi} and for
algebras in CD(3) by Kiss and Valeriote \cite{Kiss/Valeriote}. In
this paper we generalize the latter result to algebras in CD(4).

\section{Preliminaries}\label{sect:prelim}

\subsection{Constraint Satisfaction Problems and Bounded width}

Most of the terminology introduced in this section is fairly standard.
A {\em vocabulary} is a finite set of relation symbols or predicates.
In what follows, $\vocab$ always denotes a vocabulary. Every relation
symbol $R$ in $\vocab$ has an arity $r\geq 0$ associated to it. We also
say that $R$ is an $r$-ary relation symbol.

A $\vocab$-structure $\aest$ consists of a set $A$, called the {\em universe}
of $\aest$, and relations $R^{\aest}\subseteq A^r$ for every relation symbol $R\in\vocab$
where $r$ is the arity of $R$. All structures in this paper are assumed to
be {\em finite}, i.e., structures with a finite universe. Throughout the paper we
use the same boldface and slanted capital letters to denote a structure and its
universe, respectively.

A homomorphism from a $\vocab$-structure $\aest$ to a $\vocab$-structure $\best$
is a mapping $h:A\rightarrow B$ such that for every $r$-ary $R\in\vocab$ and every
$(a_1,\dots,a_r)\in R^{\aest}$, we have $(h(a_1),\dots,h(a_r))\in R^{\best}$. We will write
$\aest\rightarrow\best$, meaning that there exists a homomorphism from $\aest$
to $\best$ and say $\aest$ {\em is homomorphic to} $\best$.

\begin{definition}(Constraint Satisfaction Problems)
Let $\best$ be a finite relational structure. $\csp(\best)$ is defined
to be the set of all structures $\aest$ such that $\aest\rightarrow\best$.
Alternatively, we view $\csp(\best)$ as the computational problem asking
to decide whether a given $\vocab$-structure $\aest$ (the input) is
homomorphic to $\best$.
\end{definition}

The notion of bounded width has several equivalent formulations. In this
paper we shall base our definition on a variant of the existential $k$-pebble
game~\cite{Kolaitis/Vardi} due to Feder and Vardi~\cite{Feder/Vardi}.

\begin{definition}
Let $0\leq j<k$ be integers and let $\aest$ and $\best$ be relational structures.
A partial homomorphism from $\aest$ to $\best$ is any mapping from some subset of
the universe of $\aest$ to $\best$ that preserves all tuples of $\aest$ entirely
contained in its domain. Given two mappings $f,g$ we say that $g$ extends $f$,
denoted by $f\subseteq g$ if the domain of $f$ is a subset of that of $g$ and both
coincide over the domain of $f$.
A {\em winning strategy for the duplicator in the existential $(j,k)$-pebble game on
$\aest$ and $\best$} - or a {\em $(j,k)$-strategy} or even just a {\em strategy} if the rest
is clear - is a nonempty set $H$ of partial homomorphisms satisfying the following
two conditions:
\begin{itemize}
\item {\em Closure under subfunctions}. If $g\in H$ and $f\subseteq g$ then $f\in H$
\item {\em $(j,k)$-forth property}. If $I\subseteq J\subseteq A$ with
$|I|\leq j$ and $|J|\leq k$ and $f\in H$ with domain $I$, then there
exists $g\in H$ with domain $J$ such that $f\subseteq g$.
\end{itemize}
\end{definition}

There is a standard procedure~\cite{Larose/Zadori}, called $(j,k)$-consistency (in \cite{Larose/Zadori} they called it the (j,k)-algorithm), that given two relational structures $\aest$ and $\best$ returns, if
it exists, a $(j,k)$-winning strategy. The $(j,k)$-consistency
algorithm starts by throwing initially in $H$ all partial
homomorphisms with domain size $\leq k$. Once this is done the
procedure removes all those mappings that falsify one of the two
conditions that define a winning strategy. At the end of this
iterative process we either get a winning strategy or an empty set,
implying that such a strategy does not exist. It is not difficult to
verify that this process runs in time exponential on $k$, but
polynomial if $k$ is fixed.

Observe that every satisfiable instance has a winning strategy that consists of all $k$-ary subfunctions of the
solution.  The converse is not true. A structure $\best$ has width $(j,k)$ if the opposite always holds. More
formally,

\begin{definition}
A $\sigma$-structure $\best$ has {\em width $(j,k)$} if for every $\sigma$-structure $\aest$, if there
exists a winning $(j,k)$-strategy then $\aest$ is homomorphic to $\best$. Furthermore, $\best$ is said to
be of {\em width $j$} if it has width $(j,k)$ for some $k$ and to be of {\em bounded width} if it has
width $j$ for some $j$.
\end{definition}

It is a major open problem in the area to characterize all structures with bounded width. Up to the present
moment only width $1$ has been characterized~\cite{Feder/Vardi} (see also \cite{Dalmau/Pearson}).
Another very important question is the existence of an infinite hierarchy, i.e., whether for any $j$ there are structures of bounded width but that do not have width $j$. This is known to be true for $j=1$ but open even for $j=2$.


\subsection{Congruence distributive algebras}

In this subsection we are going to define the algebraic notions used in the paper. We assume that the reader is familiar with basic notions and results of 
universal algebra, such as algebras, varieties, congruences, clones, and so on. Good textbooks are \cite{burris/sank} and \cite{alvin}. We do not use the results of the Tame Congruence Theory (see \cite{hobby-mck}) in this paper, except for its mention in the Introduction, but it is fair to say that iteration of J\'onsson terms at the beginning of Section \ref{sec:technicalsection} was in part inspired by its basic methods. Contrary to the standard notation, we are using $\algA$, $\algB$ etc. to denote algebras, as we reserved the boldfaced letters for relational structures.

It is well known that congruences of any algebra form an algebraic lattice under the inclusion order. A famous result by B.~J\'onsson \cite{jonsson} states that an algebra $\algA$ lies in a congruence-distributive variety (congruence lattice of any algebra in the variety is distributive) if and only if there exists some $n>0$ such that $\algA$ has ternary term operations $p_0,p_1,\dots,p_n$ and the following equations hold in $\algA$:
$$
\begin{array}{ccl}
p_0(x,y,z) & \approx & x\\
p_0(x,y,z) & \approx & z\\
p_i(x,y,x) & \approx & x\\
p_i(x,x,y) & \approx & p_{i+1}(x,x,y)\hspace{1cm}\text{for all even $i$}\\
p_i(x,y,y) & \approx & p_{i+1}(x,y,y)\hspace{1cm}\text{for all odd $i$}
\end{array}
$$

We will say that an algebra ${\mathcal A}$ lies in $CD(n)$ if its fundamental operations are precisely $p_1,p_2,\dots,p_{n-1}$.

We define two operators which provide a Galois connection between algebras and relational structures: Let $\aest$ be a relational structure. Then $\Pol(\aest)$ is the clone of all operations on the universe of $\aest$ which preserve each of the relations in $\aest$. Let $\aest$ be a relational structure. Then $\Pol(\aest)$ is the clone of all operations on the universe of $\aest$ which preserve each of the relations in $\aest$. Let $\algA$ be an algebra. Then $\Inv(\algA)$ is the relational clone (set of relations closed under constructions via primitive positive formulas) of all relations on the universe of $\algA$ which are preserved by each of the operations of $\algA$. In other words, the relations of $\Inv(\algA)$ are all subuniverses of all finite powers of $\algA$.

The well-known result of Bulatov, Jeavons and Krokhin \cite{Bulatov/Jeavons/Krokhin} states that, when $\best_1$ is a relational structure all of whose relations lie in $\Inv(\Pol(\best))$, -or even in $\Inv(\Pol_{id}(\best))$, where $\Pol_{id}(\best)$ denotes the idempotent subclone of $\Pol(\best)$- then the problem $\csp(\best_1)$ is not harder than $\csp(\best)$. Therefore, we may say that an algebra $\algB$ is tractable, meaning that the problem $\csp(\best)$ is tractable for any relational structure $\best$ with relations in $\Inv(\algB)$. Having in mind the reduction to idempotent subclone, in this paper we are interested only in finite idempotent algebras and the varieties they generate, that is finite algebras in which each fundamental operation $f$ satisfies the identity $f(x,x,\dots,x)\approx x$.

\section{Main Theorem}

An algebra $\algB$ has bounded width if every structure $\best$ with relations in $\Inv(\algB)$ has bounded width. We are now ready to state the main result of this paper:

\begin{theorem}
\label{the:main}
Every algebra in CD(4) has bounded width.
\end{theorem}

Clearly, our result proves also that every algebra which has non-trivial J\'onsson terms $p_1$, $p_2$ and $p_3$ has bounded width, as adding operations to an algebra reduces the set of compatible relations, making the set of possible inputs of the related constraint satisfaction problem smaller. Therefore, in most papers in the area, when an algebra satisfies a Mal'cev-style condition, we immediately assume that the term(s) guaranteed by this condition are all the terms of the algebra.

The proof of Theorem \ref{the:main} spans the next two sections.

\section{The structure of relations}
\label{sec:technicalsection}

Recall that any algebra $\algB$ in $CD(4)$ has three basic
term operations $p_1,p_2,p_3$. We shall denote $p_2(y,x,x)$ by
$l(x,y)$ and $p_2(x,x,y)$ by $r(x,y)$. Note that the J\'onsson
equations imply that $p_1(y,x,x)=l(x,y)$ and $p_3(x,x,y)=r(x,y)$.

\begin{definition}
Let $\algD$ be a member of $CD(4)$ and $C$ a nonempty subuniverse of $\algD$. $C$ is an $l$-ideal of $\algD$ if
for every $x,y\in D$, $l(x,y)\in C$ whenever $x\in C$. The concept
of $r$-ideal is defined similarly. $\algD$ is said to be $l$-ideal free
if its only $l$-ideal is itself. $\algD$ is said so be ideal free it it
is $l$-ideal free and $r$-ideal free. The $l$-ideal of $\algD$ generated by an element $a\in D$ is the smallest $l$-ideal of $\algD$ containing $a$. An element $a\in D$ generates a minimal $l$-ideal if the generated ideal contains no proper subuniverses which are $l$-ideals of $\algD$. Analogous notions can be defined for $r$-ideals.
\end{definition}

Let $\algB$ be a finite algebra in $CD(4)$ and let $p_1,p_2,p_3$ be the J\'onsson terms of $\algB$. We shall do some preprocessing over these operations. In particular we want to guarantee that $p_1$, $p_2$, and $p_3$ besides obeying the J\'onsson identities satisfy a couple of equations more. More precisely, we need that $l(x,l(x,y))=l(x,y)$ and $r(x,r(x,y))=r(x,y)$.

We shall see how to obtain from $p_1$, $p_2$, and $p_3$, a new family of terms $p'_1$, $p'_2$, and $p'_3$ that
satisfies all required identities.

For every $x$ consider the function $l_x$ that maps every element $y$ to
$l(x,y)$. There exists some natural $n_x$ such that composing
$l_x$ with itself $n_x$ times we obtain a retraction.

We define inductively the sequence of operations
$q_1^i(x,y,z)$, $i\geq 0$ with rules: (i) $q_1^0=x$ and (ii)
$q_1^{i+1}=q_1^i(q_1(x,y,z),y,z)$. It is easy to verify by induction
that for every $i$, $q_1^i$ satisfies the identities: $q_1^i(x,x,y)=q_1^i(x,y,x)=x$, and
$q_1^i(y,x,x)=(l_x)^i(y)$. Let us fix $p'_1$ to be
$q_1^{n_1}$ with $n_1=\prod_{x\in A} n_x$. Similarly define $p'_3$ and $n_3$.
Finally, define $p'_2(x,y,z)$ as
$p_2(q_1^{n_1-1}(x,y,z),y,q_3^{n_3-1}(x,y,z))$. It is easy to verify
that $p'_1$, $p'_2$ and $p'_3$ satisfy the required identities.


From now on we fix the finite algebra $\algB$ in $CD(4)$ and the variety ${\mathcal V}={\mathcal V}(\algB)$. We stipulate that the J\'onsson terms $p_1$, $p_2$ and $p_3$ satisfy the additional equations $l(x,l(x,y))=l(x,y)$ and $r(x,r(x,y))=r(x,y)$ in ${\mathcal V}$. The following observation is going to be used a good number of times.

\begin{observation}
Let $\algB_1$ be a finite algebra in ${\mathcal V}$ and let $X$ be a subuniverse of $\algB_1$.
If $X$ is not an $l$-ideal of $\algB_1$, we can always find some
$x$ in $X$ and some $x'$ in $B_1\setminus X$ such that $l(x,x')=x'$. Same applies to $r$-ideals.
\end{observation}

\begin{lemma}
\label{prop:1} Let $\algB_1$ and $\algB_2$ be finite algebras in ${\mathcal V}$, 
$\algB_1$ $l$-ideal free, $\algD$ be a minimal $r$-ideal of $\algB_2$, and let $\algR\leq\algB_1\times\algD$ be subdirect.
If $B_1\times\{d\}\subseteq R$ for some $d\in D$ then $R = B_1\times D$. The same statement holds with $l$ and $r$ changing places.
\end{lemma}

\proof
Put $E = \{\, e\in D : B_1\times\{e\} \subseteq R \,\}$. By our assumption $E$ contains $d$, and our goal is to show that $E$ is an $r$-ideal of 
$\algB$. 
Clearly $E$ is a subalgebra. Suppose that $E$ is not an $r$-ideal.
Then there exists $e\in E$ and $e'\in B_2\setminus E$ such that $p_2(e,e,e')=e'$.
Since $e\in D$ and $D$ is an $r$-ideal of $\algB_2$, we get that $e'\in D$.
Put $C = \{\, c\in B_1 : (c,e') \in R \,\}$. As $R$ is subdirect and $e'\not\in E$, $C\neq\emptyset$ is a proper subuniverse of $B_1$. We show that $C$ is an $l$-ideal of $\algB_1$. Take $c\in C$ and $a\in B_1$. Then $(a,e), (c,e), (c,e')\in R$ and therefore $(p_2(a,c,c),p_2(e,e,e')) = (l(c,a),e')\in R$, that is $l(c,a)\in C$. This proves that $C = B_1$ which is a contradiction.
\endproof

In this section, $\algB_1$ and $\algB_2$ will always denote finite
algebras in ${\mathcal V}$ and $\algR$ will be a subdirect product
of $\algB_1$ and $\algB_2$. We shall define $G_1$, as the
subuniverse of $\algB_1^2$ that contains all tuples $(a,a')\in
B_1^2$ such that there exists some $b$ such that both $(a,b)$ and
$(a',b)$ are in $R$. We shall regard $G_1$ as a reflexive graph.

\begin{lemma}
\label{prop:2}
If $\algB_1$ is simple and $R$ is not the graph of a homomorphism $\algB_2\rightarrow \algB_1$ then
$G_1$ is connected.
\end{lemma}

\proof
Indeed, by composing $G_1$ with itself a large enough number of times we obtain a graph $G_1^*$
that has an edge precisely in those elements that are connected in $G_1$. $G_1^*$ is
a congruence and hence trivial. If $G_1^*$ is the identity then we can guarantee that $R$ is a homomorphism
$\algB_2\rightarrow \algB_1$, which is impossible. Hence we can conclude that $G_1^*$ is $B_1\times B_1$
and hence $G_1$ is connected.
\endproof

Let $X\subseteq B_1$ and $Y\subseteq B_2$, we shall say that $X$
{\em sees} $Y$ if for every $y\in Y$ there exists a tuple
$(x,y)\in R$ with $x\in X$ and similarly that $Y$ sees $X$ if for
every $x\in X$ there exists some tuple
$(x,y)\in R$ with $y\in Y$. We shall also say that
$a$ sees a set $Y$ meaning that $\{a\}$ sees $Y$.

An element $a$ of $B_1$ is said to be $2$-fan if it can see two
different elements of $B_2$. Similarly we define $2$-fan elements of
$B_2$. Obviously, $R$ is not the graph of homomorphism from $\algB_2$ to $\algB_1$ iff $B_2$ contains a 2-fan element.

\begin{lemma}
\label{l:subdirect-ideal}
Let $\algB_1$ be simple, $\algR\leq\algB_1\times\algB_2$ be subdirect,
and $\algS$ be an $r$-ideal of $\algR$. 
Assume that $\algR$ is not the graph of a homomorphism from $\algB_2$ onto $\algB_1$, and that $\algS$ is the graph of a homomorphism of $\algD = \pi_2(\algS)$
onto $\algC = \pi_1(\algS)$. 
Then $r(c,a) = c$ for all $c\in C$ and $a\in B_1$. The analogous statement with $l$ replacing $r$ everywhere also holds.
\end{lemma}

\proof Since $\algR$ is not the graph of a homomorphism and $\algB_1$ is simple, $G_1$ is connected. Clearly, $r(c,c) = c$ for any $c\in C$. 
By using the connectivity of $G_1$ it is enough to show that $r(c,a') = c$
whenever $r(c,a) = c$ and $(a,a')\in G_1$.
As $c\in C$, there exists $d\in D$ so that $(c,d)\in S$.
Let $b\in B_2$ be such that $(a,b),(a',b)\in R$.
Then $r((c,d),(a',b)) = (r(c,a'), r(d,b)) \in S$
and $r((c,d),(a,b)) = (r(c,a),r(d,b)) = (c,r(d,b)) \in S$.
Since $S$ is the graph of homomorphism we get that $r(c,a') = c$.
\endproof

\begin{lemma}
\label{le:3} Let $\algB_1$ and $\algB_2$ be finite algebras in
${\mathcal V}$, where $\algB_1$ is simple and ideal free, and let
$\algR$ be a subdirect product of $\algB_1$ and $\algB_2$. If
$B_2$ contains a $2$-fan element then it also contains an element
that sees the whole of $B_1$.
\end{lemma}

\proof Let us assume that $\algR$ is a counterexample to the
statement with $|B_1|+|B_2|$ as small as possible. We shall do a
separate analysis depending on whether $B_2$ has a proper ideal and/or congruence. In all cases we shall reach a contradiction.

CASE 1: $\algB_2$ is not ideal free. Let $Y$ be a proper $l$-ideal of $\algB_2$. We first prove that $Y$ sees the whole of $B_1$ by showing that the
subset $Z$ of $B_1$ that contains all elements seen by $Y$ is an
$l$-ideal of $B_1$. Indeed, let $a,a'$ be any elements of $B_1$ with
$a\in Z$. Let $b,b'$ be elements of $B_2$ seen by $a,a'$
respectively. Since $a\in Z$ we can assume that $b\in Y$. By
applying $l$ to $(a,b)$ and $(a',b')$ we obtain $(l(a,a'),l(b,b'))$.
Since $l(b,b')\in Y$ we conclude that $l(a,a')\in Z$. 

So, the projection $S$ of $R$ to $B_1\times Y$ is subdirect and a proper $l$-ideal of $\algR$. As $|B_1|>1$, the ideal freeness of $\algB_1$ implies that there exist elements $a,b\in B_1$ such that $l(a,b)\neq a$. According to Lemma~\ref{l:subdirect-ideal}, this means that $S$ contains an element which sees two elements of $B_1$, a contradiction with minimality of $(B_1,B_2,R)$. Analogously we prove that $B_2$ can have no $r$-ideals.

CASE 2: $\algB_2$ is not simple. Let $\theta$ be a non trivial
congruence of $B_2$. Consider now the relation $S$ defined as
$$\{(a,b/\theta) \; | \; (a,b)\in R\}$$ If $b$ is any element in
$B_2$ that has $2$-fan in $R$, then $b/\theta$ has $2$-fan in $S$.
By the minimality of $R$ and Lemma \ref{prop:1},
$S=B_1\times B_2/\theta$. Let $Y$ be
$b/\theta$ (now regarded as a subset of $B_2$). We have that $Y$
sees the whole of $B_1$ and contains a $2$-fan element, namely $b$.
By the minimality of $R$, $B_2$ contains an element that sees the
whole of $B_1$, a contradiction.

CASE 3: $\algB_2$ is ideal free and simple. Select any $2$-fan
element $b$ of $B_2$ and construct the sequence $Y^0,Y^1,\dots$,
defined inductively by the following rules: (i) $Y^0=\{b\}$, and
(ii) $Y^{i+1}$ is the set of elements seen by $Y^i$. By Lemma
\ref{prop:2} $G_1$ is connected and therefore the sequence reaches
at some point $B_1$ (and hence also $B_2$). Let $Y^i$ be the latest
element of the sequence before any of $B_1$ or $B_2$ occurs. We
consider two possibilities: If $Y^i\subseteq B_2$, then it must
contain $b$. Consider the relation $S$ defined as $R\cap (B_1\times
Y^i)$. By the minimality of $R$, the relation $S$ must contain some
element that sees the whole of $B_1$. The very same element should
also see the whole of $B_1$ in $R$, a contradiction. If
$Y^i\subseteq B_1$ we see that if $Y^i$ has no 2-fan elements, then
$Y^{i-1}=Y^{i+1}=B_2$, which contradicts the choice of $Y^i$. Hence
we conclude that there exists an element of $B_1$ that sees the
whole of $B_2$. By Lemma \ref{prop:1} $R=B_1\times B_2$ and we
are done.
\endproof

The next lemma is quite trivial, but its use is going to be essential
in the construction of substrategies.

\begin{lemma}
\label{l:min-ideal}
Let $\algB_1$ and $\algB_2$ be algebras in ${\mathcal V}$ and let $\algR$ be a subdirect product of $\algB_1$ and $\algB_2$.
\begin{enumerate}
\item If $(a,b)\in R$ generates a minimal $l$-ideal in $\algR$, then $a$ generates a minimal $l$-ideal in $\algB_1$ and $b$ generates a minimal $l$-ideal in $\algB_2$.
\item If $a\in B_1$ generates a minimal $l$-ideal in $\algB_1$, then
there exists $b\in B_2$ such that $(a,b)\in R$ and $(a,b)$ generates
a minimal $l$-ideal in $\algR$.
\end{enumerate}
Same statements holds for $r$-ideals.
\end{lemma}

\proof
To prove statement~$(1)$, assume that $(a,b)\in R$ generates a minimal $l$-ideal in $\algR$. Let $\algC$ be the ideal of $\algB_1$ generated by $a$. If $\algC$ is not minimal, then there exists an element $c\in C$ that does not generate $a$. However, as $\algR$ is subdirect, there is $d\in B_2$ so that $(c,d)\in R$ is generated by $(a,b)$ (just follow the steps in generation of $c$ by $a$, start from $(a,b)$ and whenever a constant $x\in B_1$ is used, replace it by $(x,y)\in R$). Therefore $(c,d)$ cannot generate $(a,b)$ in $\algR$, which is a contradiction with the choice of $(a,b)$.

To prove statement~$(2)$, assume that $a\in\algB_1$ generates a minimal ideal $C$ of $\algB_1$. 
Since $\algR$ is subdirect, there exists $b\in B_2$ such that $(a,b)\in R$.
Let $\algS$ be the ideal of $\algR$ generated by $(a,b)$, and choose 
$(c,d)\in S$ that generates a minimal ideal $\algS'$ of $\algR$ (any element of a minimal ideal $S'\subseteq S$ generates it, of course). Since $(c,d)\in S$, then $c\in C$. Then since $C$ is minimal, we see that $C$ is also the ideal of $\algB_1$ generated by $c$. So, as $a$ is in the ideal generated by $c$, we can generate an element $(a,b')\in S'$. This element generates $S'$, a minimal ideal in $\algR$.
\endproof

\begin{lemma}
\label{l:subdirect-ideal-product}
Let $\algB_1$ be simple, ideal free
member of ${\mathcal V}$ and let $\algB_2\in{\mathcal V}$. Let $\algR\leq\algB_1\times\algB_2$ be a subdirect 
product that is not the graph of a homomorphism from $\algB_2$ onto $\algB_1$. 
If $\algS$ is a minimal $l$-ideal of $\algR$ or a minimal $r$-ideal of $\algR$, then 
$S = B_1\times\pi_2(\algS)$.
\end{lemma}

\proof
Put $\algC = \pi_1(\algS)$ and $\algD = \pi_2(\algS)$.
By Lemma~\ref{l:min-ideal}, $\algC$ and $\algD$ are minimal $l$-ideals
of $\algB_1$ and $\algB_2$, respectively. 
However, $\algB_1$ is $l$-ideal free, thus $C = B_1$. 
We cannot have $p_2(c,c,a) = c$ for all $a,c\in B_1$ because then
every one-element subset of $\algB_1$ would be an $l$-ideal of $\algB_1$.
Therefore, by Lemma~\ref{l:subdirect-ideal}, $\algS$ is not the graph of a homomorphism. 
Now $\algS\leq B_1\times D$ is subdirect, then by Lemma~\ref{le:3}, 
there exists $d\in D$ such that $B_1\times\{d\}\subseteq S$. 
Now using Lemma~\ref{prop:1} we get that $S = B_1\times D$.
\endproof

\section{Proof of Theorem~\ref{the:main}}

If $H$ is a strategy and $I=\{a_1,\dots,a_i\}$ is a subset of $A$ with $i\leq k$ we denote
by $H_I$ the subset of $H$ that contains precisely all the mappings with domain $I$. An alternative
but essentially equivalent view is to fix an order on the elements of $I$, say $a_1<a_2<\cdots<a_i$,
and to regard $H_I$, or rather $H_{a_1,\dots,a_i}$, as the $i$-ary relation on $B$
$$\{(f(a_1),\dots,f(a_i)) \; | \; f\in H, \dom(f)=I\}$$

In what follows we shall assume that every relation of $\best$ is in $\Inv(\algB)$. In this case, it is easy to verify -and widely known- that the strategy $H$
returned by the $(j,k)$-consistency procedure satisfies the following property: for every
$a_1,\dots,a_i$, $H_{a_1,\dots,a_i}$ is in $\Inv(\algB)$. We shall say, somehow abusing notation,
that $H\in\Inv(\algB)$. We shall apply some transformations to the
winning strategy returned by the $(j,k)$-consistency algorithm. In all our transformations this property will be maintained.

We shall also fix the values of $j$ and $k$. From now on $j$ is assumed to be $k-1$ and $k$ is the maximum
between $3$ and the largest of the arities of signature $\sigma$. Observe that if $j=k-1$ then the
$(j,k)$-forth property can be rephrased as follows:
If $f\in H$ with domain $|I|<k$ and $a\in A$ then there exists some $g\in H$ defined on $I\cup\{a\}$
such that $f\subseteq g$. We will call the $(k-1,k)$-forth property as the $k$-forth property from now on.

In order to simplify notation we
shall omit the parameter $j$ and we shall speak of $k$-winning strategy, $k$-consistency algorithm and so on.

\begin{lemma}
\label{le:idealreduction}
Let $\aest$ and $\best$ be $\sigma$-structures such that $\best\in\Inv(\algB)$ and let $H\in\Inv(\algB)$ be a
$k$-strategy. Then there exists a $k$-strategy $H'\in\Inv(\algB)$ where for every $a\in A$, $H'_a$ is ideal free.
\end{lemma}
\proof
This proof shamelessly duplicates that of Lemma 3.14 in~\cite{Kiss/Valeriote}. We include it here for the
sake of completeness. It must be pointed out that the consideration of strategies, instead of
relational width, as in~\cite{Kiss/Valeriote} simplifies slightly the job. Indeed, in~\cite{Kiss/Valeriote}
Lemma 3.14 has be complemented with Lemmas 3.16 and Corollaries 3.15 and 3.17 in order to achieve the same result.

For simplicity of notation we shall use integers to denote the elements of $A$.
Let us assume that $H_1$ has a proper ideal $X$.
We shall obtain a new $k$-strategy $H'$ such that
$H'_1=X$ in the following way:
\begin{itemize}
\item In the first stage we place
in $H'$ every mapping $g\in H$ with $1$ in its domain such that $g(1)\in X$ and
every one of its subfunctions $f\subseteq g$.
\item In the second stage we include in $H'$ all the mappings $f$ of $H$ such that the domain of $f$ has exactly $k$ elements and
every one of its proper subfunctions was included in the first stage.
\end{itemize}

It is routine to verify that $H'$ is nonempty, closed under subfunctions and in $\Inv(\algB)$.
It has to be proved that $H'$ has the $k$-forth property. We shall present the proof in the
case that $X$ is an $l$-ideal, which uses operation $p_1$. The proof for $X$ being an $r$-ideal
is obtained analogously by using operation $p_3$.

Let $f\in H'$ be a mapping with domain
$I$ with $|I|<k$ and let $i$ be any element of $A$. We have to prove that there exists
some mapping $g\in H'$ defined on $i$ that extends $f$.
We shall observe first that the only challenging case is when $|I|=k-1$. Indeed, if $|I|<k-1$
the extension $g$ is obtained in the following way: First observe that mapping $f$ can
only be added to $H'$ in the first stage. Hence there exists some mapping $h\in H$
with domain $\{1\}\cup I$ such that $h(1)\in X$. Hence by the $k$-forth property, there exists some
extension $g\in H$ of $h$ defined on $\{1,i\}\cup I$. Since $h(1)=g(1)\in X$, $g$ is also included in the first stage, as is its restriction to $I\cup\{i\}$, which extends $f$.

So for now we shall assume that $|I|=k-1$. The case $1\in
I\cup\{i\}$ is also straightforward. Since $f$ has to be included in
$H$ in the first stage, there exists some $g\in H$ with domain
$I\cup\{1\}$ such that $g(1)\in X$. If $i=1$ then we are done.
Otherwise, $1\in I$ and $f=g$. Mapping $f$ can be extended to
$I\cup\{i\}$. The obtained mapping belongs to $H'$ because it is
necessarily included in the first stage.

Hence, we can assume  that $1\not\in I\cup\{i\}$.
This turns out to be the more complicated case. Let us set, for ease of notation, that
$I=\{2,3,4,\dots,k\}$ and that $i=k+1$. We shall show that there exists some
$b_{k+1}$ such that $(b_2,\dots,b_{k+1})\in H'_{2,3,\dots,k+1}$.

By the $k$-forth property of $H$ there exists some extension
$(b_2,b_3,b_4\dots,b_{k},u_{k+1})$ $\in H$ of $f$. Also the mapping
$f$ has been included in the first stage due to the fact that there
exists some extension $g$ of $f$ with $g(1)=b_1\in X$. We conclude
that $(b_1,b_2,b_3,b_4\dots,b_k)\in H_{1,\dots,k}$. By applying
successively the closure under subfunctions and the $k$-forth
property of $H$ we conclude that $H$ contains some assignments
$(b_1,b_3,b_4\dots,b_k,v_{k+1})$, $(v_2,b_3,b_4\dots,b_k,v_{k+1})$,
$(b_1,b_2,b_4,\dots,b_k,w_{k+1})$, and
$(b_2,w_3,b_4,\dots,b_{k},w_{k+1})$ (the domains of the mappings are
implicitly indicated by the indexes). Let $b_{k+1}$ be
$p_1(u_{k+1},v_{k+1},w_{k+1})$. By applying $p_1$ to
$(b_2,b_3,b_4\dots,b_{k},u_{k+1})$, $(v_2,b_3,b_4\dots,b_{k},$
$v_{k+1})$, and $(b_2,w_3,b_4\dots,b_{k},w_{k+1})$ we conclude that
the tuple $(b_2,b_3,b_4\dots,b_{k},b_{k+1})$ belongs to $H$. We need
to show that for all $2\leq i\leq k$,
$(b_2,\dots,b_{i-1},b_{i+1},\dots,b_{k+1})$ was included in the
first stage, or equivalently, that there exists some $c_1\in X$ such
that $(c_1,b_2,\dots,b_{i-1},b_{i+1},\dots,b_{k+1})\in H$. There are
a number of cases to consider:
\begin{itemize}
\item If $i=k+1$ then the tuple $(b_2,\dots,b_{k})$ extends to $(b_1,\dots,b_{k})$, as required.
\item If $i=2$ then by the properties of $H$ we can conclude that $H$ contains some tuples
$(x_1,b_3,\dots,b_k,u_{k+1})$ and $(b_1,y_3,b_4,\dots,b_{k},w_{k+1})$. Applying $p_1$ to these tuples
along with the tuple $(b_1,b_3,b_4,\dots,v_{k+1})$ we obtain the tuple $(l(b_1,x_1),b_2,\dots,b_k)$.
Since $X$ is an $l$-ideal, $l(b_1,x_1)\in X$.
\item If $i=3$ or $3<i<k+1$ then small variations of the previous argument will work.
\end{itemize}

By repeated application of the procedure we shall obtain the required strategy.
\endproof



\begin{lemma}
\label{l:simple-platoo}
Let $\aest$ and $\best$ be $\sigma$-structures such that $\best\in\Inv(\algB)$ and let $H\in\Inv(\algB)$ be a $k$-strategy.
Assume that for all $i\in A$, $|H_i|\geq 2$ and $\algH_i$ is ideal free.
Then there exists a nonempty subset $M\subseteq A$ and 
maximal congruences $\vartheta_m$ of $\algH_m$ for all $m\in M$ that satisfy the following statements.
\begin{enumerate}
\item For any pair $m_1,m_2\in M$ of distinct elements, $\algH_{m_1,m_2}/(\vartheta_{m_1}\times\vartheta_{m_2})$ is the graph of an isomorphism $\tau_{m_1,m_2} : \algH_{m_1}/\vartheta_{m_1} \to \algH_{m_2}/\vartheta_{m_2}$.
\item For any $m\in M$ and $n\in A\setminus M$,
$\algH_{n,m}/(0_{\algH_n}\times\vartheta_m) = \algH_n\times\algH_m/\vartheta_m$.
\item By setting $\tau_{m,m}$ to be the identity on $\algH_m/\vartheta_m$ for all $m\in M$, then $\tau_{m_1,m_2}\circ\tau_{m_2,m_3} = \tau_{m_1,m_3}$ for all $m_1,m_2,m_3\in M$. 
\end{enumerate}
\end{lemma}

\proof
Let $M\subseteq A$ be of maximal size with respect to satisfying statement~$(1)$. Then $M$ is nonempty, as any one-element subset of $A$ satisfies that condition.
Assume that statement~$(2)$ is not satisfied, that is there exist elements $m\in M$ and $n\in A\setminus M$ and such that $\algH_{n,m}/(0_{\algH_n}\times\vartheta_m)$ is not the direct product. Then by Lemma~\ref{l:subdirect-ideal-product}, 
$\algH_{n,m}/(0_{\algH_n}\times\vartheta_m)$ is the graph of a homomorphism $\varphi : \algH_n\to\algH_m/\vartheta_m$. 
Put $\vartheta_n = \ker\varphi$ and $\tau_{n,m} = \varphi/\vartheta_n$.
Clearly, $\vartheta_n$ is a maximal congruence of $\algH_n$ and $\tau_{n,m}$ is
an isomorphism. Let $m'\in M\setminus\{m\}$ be any element. Since $H$ is also a $(2,3)$-strategy, 
$\pi_{1,2}(H_{n,m,m'}) = H_{n,m}$, $\pi_{2,3}(H_{n,m,m'}) = H_{m,m'}$ and 
$\pi_{1,3}(H_{n,m,m'}) = H_{n,m'}$.
Now the projection of $\algH_{n,m,m'}/(\vartheta_n\times\vartheta_m\times\vartheta_{m'})$ onto the first two and last two coordinates yield the graphs of the isomorphisms $\tau_{n,m}$ and $\tau_{m,m'}$, respectively, therefore the projection onto the first and last coordinate yields the graph of the isomorphism $\tau_{n,m'} = \tau_{n,m}\circ\tau_{m,m'}$. This proves that $M\cup\{n\}$ also satisfies statement~$(1)$, which is a contradiction. Note, that the last argument of the proof proves statement~$(3)$ as well.
\endproof

\begin{lemma}\label{le:simplereduction}
With the assumptions of Lemma~\ref{l:simple-platoo}, for all $m\in M$ let
$C_m$ be a congruence class of $\vartheta_m$ that correspond to each other under the isomorphisms $\tau_{m_1,m_2}$. Let $G$ be the set of all functions $g\in H$ that satisfy the following conditions
\begin{enumerate}
\item $g(m)\in C_m$ for all $m\in\dom(g)\cap M$, and
\item $g$ generates a minimal $r$-ideal of $\algH_{\dom(g)}$.
\end{enumerate}
Then $G$ is a $k$-strategy. The $k$-strategy $\overline{G}$ generated by $G$
(with functions $\overline{G}_K = \sg_{\prod_{k\in K}\algH_k}(G_K)$ for all $K\subseteq A$) is a $k$-strategy in $\Inv(\algB)$ that also satisfies condition~$(1)$.
\end{lemma}

\proof
Clearly $G$ is nonempty, as for any element $m\in M$ and element $b\in\algH_m$ generates a minimal $r$-ideal of $\algH_m$ and therefore the function $g$ with domain $\{m\}$ with $g(m) = b$ is in $G$. Clearly, the set of functions satisfying condition~$(1)$ is closed under subfunctions, and also the ones satisfying condition~$(2)$ because of Lemma~\ref{l:min-ideal}. We need to prove the $k$-forth property. 

Let $f\in G$ with $|\dom(f)|<k$ and choose $i\in A\setminus\dom(f)$.
Put $J = \dom(f)$ and $K = \dom(f)\cup\{i\}$.
Note, that $\algH_K$ is the subdirect product of $\algH_J$ and $\algH_i$,
therefore, by Lemma~\ref{l:min-ideal}, there exists a function $h\in H_K$ such that $h$ generates a minimal $r$-ideal in $\algH_K$ and $h|_J = f$.
If $i\not\in M$, then $h\in G$ and we are done. So assume that $i\in M$.

If $J\cap M\neq\emptyset$, then for $j\in J\cap M$, $(h(i),h(j))\in\algH_{i,j}$, and as $h(j)\in C_j$ we get that $h(i)\in C_i$ as desired. 
So we can assume that $J\cap M = \emptyset$.

Let $\algD = \algH_i$, $\algE = \algH_J$ and $\algR = \algH_{i,J}$,
and $\algS$ be the minimal $r$-ideal of $\algR$ generated by the
element $h$. 
Put $\hat{\algD} = \algD/\vartheta_i$, 
$\hat\vartheta =\vartheta_i\times 0_{\algH_J}$, 
$\hat{\algR} = \algR/\hat\vartheta$ and
$\hat{\algS} = \algS/\hat\vartheta$. Clearly $\hat{\algD}$ is simple
ideal free, $\hat{\algR}$ is a 
subdirect product of $\hat{\algD}$ and $\algE$, and $\hat{\algS}$ is a minimal $r$-ideal of $\hat{\algR}$.

Assume first that $\hat{\algR}$ is the graph of a homomorphism from $\algE$ onto $\hat{\algD}$ with $\theta$ as its kernel, then $\theta$ is a maximal (coatom) congruence in Con$(\algE)$, since $\hat{\algD}$ is simple. If we denote the kernels of projections as $\eta_j$, for $j\in J$, then 
$$\theta=\theta\vee\bigwedge\limits_{j\in J}\eta_j=\bigwedge\limits_{j\in J}(\theta\vee\eta_j)$$
(by the distributivity of congruences). As $\theta$ is maximal, there exists $j\in J$ so that $\theta\geq\eta_j$. Therefore, $\pr_{i,j}\hat{\algS}$, which is an $r$-ideal of $\pr_{i,j}\hat{\algR}$, is the graph of a homomorphism, while $\pr_{i,j}\hat{\algR}=\algH_{i,j}/(\vartheta_i\times 0_{\algH_j})$ is not the graph of a homomorphism, according to the properties of $M$. Therefore, Lemma~\ref{l:subdirect-ideal} would suggest that every element of $\pr_i\hat{\algS}$ is an $r$-ideal of $\hat{\algD}$. This contradiction implies that $\hat{\algR}$ is not the graph of a homomorphism. 

Using Lemma~\ref{l:subdirect-ideal-product} for $\hat{\algR}$ and $\hat{\algS}$ we get that 
$\hat S = \hat D\times \pi_2(\hat S) = \hat D\times\pi_2(S)$. 
However, $f = h|_{J}\in\pi_2(S)$. Therefore there exists an element $g\in S$
such that $g(i)\in C_i$ and $g|_J = f$. As $S$ was a minimal $r$-ideal, 
$g\in G$.

The last statement of the lemma is obvious (always works), but since $C_m$ for all $m\in M$ is a subalgebra, the generated strategy also satisfies condition~(1). 
\endproof

We are in a situation to prove Theorem~\ref{the:main}

\proof (of Theorem \ref{the:main}) 

By repeated application of Lemmas~\ref{le:idealreduction}, and~\ref{le:simplereduction} (before each application of Lemma~\ref{le:simplereduction} we project the new $H$ to coordinates which are not already singletons) we can construct
a winning strategy where for every $i\in A$, $H_i$ is a singleton. Since the arity of every relation
is at most $k$ the mapping $s$ sending $i$ to the only element in $H_i$ is a solution.
\endproof

\section{Acknowledgements}

The authors wish to thank Ralph McKenzie for many helpful suggestions.


\begin{thebibliography}{99}

\bibitem{Barto/Kozik/Niven} L.~Barto, M.~Kozik and T.~Niven, {\em CSP dichotomy holds for digraphs with no sources and no sinks (a positive answer to the conjecture of Bang-Jensen and Hell),} in preparation.\\

\bibitem{BIMMVW} J.~Berman, P.~Idziak, P.~Markovi\'c, R.~McKenzie, M.~Valeriote and R.~Willard, {\em  Varieties with few subalgebras of powers,} in preparation.\\

\bibitem{burris/sank} S.~Burris and H.~P.~Sankappanavar, {\em A course in universal algebra,} Graduate Texts in Mathematics, No. 78, Springer-Verlag, New York-Berlin, 1981.\\

\bibitem{Bulatovrelational} A.~Bulatov,
{\em A Graph of a Relational Structure and Constraint Satisfaction Problems,}
In Proceedings of the 19th IEEE Annual Symposium on Logic in Computer Science (LICS'04) Turku, Finland, 2004.\\

\bibitem{Bulatov2semilattices} A.~Bulatov, {\em Combinatorial problems raised from 2-semilattices,}
Journal of Algebra {\bf 298} no. 2 (2006), 321--339.\\

\bibitem{Bulatov/Jeavons/Krokhin} A.~Bulatov, P.~Jeavons and A.~Krokhin,
{\em Classifying the complexity of constraints using finite algebras,}
SIAM J. Comput. {\bf 34} (2005), no. 3, 720--742 (electronic).\\

\bibitem{Bulatov/Jeavons/Krokhin:survey} A.~Krokhin, A.~Bulatov and P.~Jeavons,
{\em The complexity of constraint satisfaction: an algebraic approach,}
Notes taken by Alexander Semigrodskikh. NATO Sci. Ser. II Math. Phys. Chem., {\bf 207}, Structural theory of automata, semigroups, and universal algebra, 181--213, Springer, Dordrecht, 2005.\\

\bibitem{Bulatov/Krokhin/Larose} A. Bulatov, A. Krokhin and B. Larose, {\em Dualities for constraint satisfaction problems,} submitted for publication.\\

\bibitem{Dalmau} V.~Dalmau, {\em Generalized majority-minority operations are
tractable,}  Log. Methods Comput. Sci. {\bf 2} (2006),  no. 4, 4:1, 14 pp.\\

\bibitem{Dalmau/Pearson} V.~Dalmau and J.~Pearson, {\em Closure functions and width 1 problems,}\\

\bibitem{Feder/Vardi} T.~Feder and M.~Y.~Vardi, {\em The computational structure of monotone monadic SNP and constraint satisfaction: a study through Datalog and group theory,} SIAM J. Comput. {\bf 28} (1999), no. 1, 57--104 (electronic).\\

\bibitem{hobby-mck}
D.~Hobby and R.~McKenzie, {\em The Structure of Finite Algebras,} Contemporary Mathematics Series Vol.\ 76, American
Mathematical Society, Providence, RI, 1991.\\

\bibitem{IMMVW} P.~Idziak, P.~Markovi\'c, R.~McKenzie, M.~Valeriote and R.~Willard, {\em Tractability and learnability arising from algebras with few subpowers,} LICS 2007.\\

\bibitem{Jeavonsetal} P.~Jeavons, D.~Cohen, and M.~Gyssens, {\em
Closure properties of constraints, } J. ACM. {\bf 44} no. 4 (1997),
527--548. \\

\bibitem{jonsson} B.~J\'{o}nsson, {\it Algebras whose congruence lattices are distributive,} Math. Scand. {\bf 21} (1967), 110--121.\\

\bibitem{Kiss/Valeriote} E.~Kiss and M.~Valeriote, {\em On tractability and congruence distributivity,} Logical Methods in Computer Science, Vol. {\bf 3} (2:6) 2007, 20 pages.\\

\bibitem{Kolaitis/Vardi} P.~G.~Kolaitis and M.~Y.~Vardi, {\em A game-theoretic approach to constraint satisfaction,} AAAI-2000/IAAI-2000 Proceedings (Austin, TX), 175--181, MIT Press, Cambridge, MA, 2000.\\

\bibitem{Larose/Tesson} B.~Larose and P.~Tesson, {\em Universal Algebra and Hardness Results for Constraint Satisfaction Problems,} ICALP 2007.\\

\bibitem{Larose/Zadori} {\em Bounded width problems and algebras,} Algebra Universalis {\bf 56} no. 3-4 (2007), 439--466.\\

\bibitem{mar-mck} M.~Mar\'oti and R.~McKenzie, {\em Existence theorems for weakly symmetric operations,} to appear in Algebra Universalis.\\

\bibitem{alvin} R.~McKenzie, G.~McNulty and W.~Taylor, {\em Algebras, lattices, varieties. Vol. I,} The Wadsworth \& Brooks/Cole Mathematics Series, Wadsworth \& Brooks/Cole Advanced Books \& Software, Monterey, CA, 1987.\\

\end{thebibliography}
\end{document}